\font\script=eusm10.
\font\sets=msbm10.
\font\stampatello=cmcsc10.
\font\symbols=msam10.
\font\piccolo=cmr8.
\font\piccolissimo=cmr5.

\def\qc{{q_{_c}}}
\def\qone{{q_{_1}}}
\def\qtwo{{q_{_2}}}
\def\qoneinv{{\overline{q}_{_1}}}
\def\sgn{\hbox{\rm sgn}}

\def\defineq{\buildrel{def}\over{=}}

\def\square{\hbox{\vrule\vbox{\hrule\phantom{s}\hrule}\vrule}}
\def\1{\hbox{\bf 1}}
\def\C{\hbox{\sets C}}
\def\N{\hbox{\sets N}}

\def\R{\hbox{\sets R}}
\def\Z{\hbox{\sets Z}}
\def\doublesum{\mathop{\sum \sum}}

\def\integrale{\mathop{\int}}

\def\Corr{\hbox{\script C}}

\def\EssBdd{\hbox{\symbols n}\,}

\par
\centerline{\bf ON THE CORRELATIONS, SELBERG INTEGRAL AND SYMMETRY}
\centerline{\bf OF SIEVE FUNCTIONS IN SHORT INTERVALS, III}
\bigskip
\centerline{by G.Coppola\footnote{$^1$}
{\piccolo titolare di un Assegno \lq \lq Ing.Giorgio Schirillo\rq \rq \thinspace dell'Istituto Nazionale di Alta Matematica (Fellow \lq \lq Ing.Giorgio Schirillo\rq \rq \thinspace of the Istituto Nazionale di Alta Matematica).} and M.Laporta}
\bigskip
{
\font\eightrm=cmr8
\eightrm {
\par

{\bf Abstract.} An arithmetic function $f$ is called a sieve function of range $Q$, if it is
the convolution product of the constantly $1$ function and $g$  such that 
$g(q)\ll_{\varepsilon} q^{\varepsilon}$, $\forall\varepsilon>0$, for $q\leq Q$, and
$g(q)=0$ for $q>Q$. Here
we establish a new result on the autocorrelation of $f$ 
by using a famous theorem on bilinear forms of Kloosterman fractions by Duke, Friedlander and Iwaniec. In particular, for such correlations we obtain  non-trivial
asymptotic formul\ae\ that are actually unreachable by the standard approach 
of the distribution of $f$ in the arithmetic progressions. Moreover, we apply our asymptotic formul\ae\ to obtain new bounds for the so-called Selberg integral and symmetry integral of $f$, which are basic tools for the study of the distribution of $f$ in short intervals.
}
\footnote{}{\par \noindent {\it Mathematics Subject Classification} (2000): 11N37, 11N25.}
}
\bigskip
\par
\noindent {\bf 1. Introduction and statement of the results}
\bigskip

A basic tool for the study of the distribution of an arithmetic function $f:\N\to\C$ in short intervals is the so-called {\it Selberg integral}  of $f$, that is
$$
J_f(N,h)\defineq \int_{N}^{2N} \Big| \sum_{x<n\le x+h}f(n) - M_f(x,h)\Big|^2\,{\rm d}x,
$$
where $M_f(x,H)$ is the (short interval) {\it mean-value} of $f$ and $h,N\in\N$ are such that $h=o(N)$, as $N\to \infty$. Indeed, non-trivial bounds for 
$J_f(N,h)$ might yield results on the distribution of $f$ in
{\it almost all} the short intervals $(x,x+h]$, i.e. for all $x\in [N,2N]\cap\N$ with $o(N)$ exceptions.
\noindent
On the other side, the symmetry properties of $f$ in almost all the short intervals are linked to the {\it symmetry integral} of $f$,
$$
I_f(N,h)\defineq \int_{N}^{2N} \Big| \sum_{|n-x|\le h}\sgn(n-x)f(n)\Big|^2\,{\rm d}x,
$$
where the {\it sign} function is defined as $\sgn(0)\defineq 0$, and $\sgn(t)\defineq |t|/t$ if $t\not=0$. 

The aim of the present paper is to continue the study of $J_f(N,h)$ and $I_f(N,h)$, started in [C1] and considered also in [C2], [CL1] and [CL2],  for a 
{\it sieve function} $f$ of {\it range} $Q\ll N$, meaning
that its {\it Eratosthenes transform} $g\defineq f\ast\mu$  is supported in $[1,Q]$ and $g$ is {\it essentially bounded}, i.e. $g(q)\ll_{\varepsilon} q^{\varepsilon}$ ($\forall\varepsilon>0$). Here, we recall that $\ll$ is Vinogradov's notation, synonimous to Landau's $O$-notation. In particular, $\ll_{\varepsilon}$ means that the implicit constant might depend on an arbitrarily small $\varepsilon>0$, which might change at each occurrence.
\noindent
Since by the M\"obius inversion formula one has
$$
f(n)=(g\ast \1)(n)=\sum_{{d|n}\atop {d\le Q}}g(d),
$$
where $\1(n)=1$ for all $n\in\N$, then $g$ is essentially bounded 
if and only if so is $f$. Moreover, since 
$$
{1\over x}\sum_{n\le x}f(n)= {1\over x}\sum_{d\le Q}g(d)\Big[ {x\over d}\Big] = \sum_{d\le Q}{{g(d)}\over d} + O\Big( {1\over x}\sum_{d\le Q}|g(d)|\Big), 
$$
where $[t]$ is the {\it integer part} of $t\in\R$ (hereafter, in sums over positive integers like $\sum_{a\le x}1$ it is implicit that $a\ge 1$), we expect the (short interval) mean-value of $f$ to be
independent of $x$, namely given by (see [CL2] for further comments)
$$
M_f(h)\defineq h\sum_{d\le Q}{{g(d)}\over d}.
$$
					%PAGE 2
Before stating our results, let us introduce some auxiliary notation and convention. When $h=[N^{\theta}]$, with $\theta\in [0,1]$, 
we refer to $\theta$ as the {\it width} of the short interval $[x-h,x+h]$ or $(x,x+h]$.
We adopt the convention that $\theta<\theta_0$ (resp. $\theta>\theta_0$) means $\theta\le\theta_0-\delta$ (resp. $\theta\ge\theta_0+\delta$) for some absolute constant $\delta>0$. Further, we say that $f$ has {\it level} $\lambda \in [0,1]$ if it is a sieve function of range $Q=[N^{\lambda}]$, and for $\lambda$ we adopt the same convention as for the width. Finally, given the arithmetic functions $\phi_1$ and $\phi_2$, we write
$\phi_1(n)\EssBdd \phi_2(n)$ to mean that $\phi_1(n)\ll_{\varepsilon} n^{\varepsilon} \phi_2(n)\ \forall \varepsilon>0$ (as $n\to \infty$).
\bigskip

\noindent {\stampatello Theorem 1.} {\it Fix  a small $\delta>0$. If $f:\N \rightarrow \R$ has level $\lambda\in(1/2,1)$, then
$$
J_f(N,h)\EssBdd Nh+N^{\delta}Q^{2-\Delta}h^2+N^{1-2\delta/3}h^2+Qh^2,
$$
$$
I_f(N,h)\EssBdd Nh+N^{\delta}Q^{2-\Delta}h^2+N^{1-2\delta/3}h^2,
$$
as $N\to \infty$, where $Q=[N^{\lambda}]$, $h=[N^{\theta}]$ with $\theta\in (0,1/2)$, and $\Delta=1/48$.
} 
\bigskip

\noindent
{\stampatello Remark.} In [C1] the above inequalities hold 
with $\Delta=0$ (for a small $h$). In particular, such inequalities yield the  non-trivial bound
$N^{1-\varepsilon}h^2$ for both integrals 
$J_f(N,h)$ and $I_f(N,h)$ with $f$ of level $\lambda <(1+\theta)/2$ and for any width $\theta\in(0,1)$ (see Corollary 1.1 of [C1]), whereas the previous theorem 
holds for $\lambda>1/2$. By combining Theorem 1 above with the results given by Corollary 1.1 of [C1] we immediately obtain the following non-trivial bounds that however improve on [C1] estimates only in very short intervals, namely $\theta\in(0,1/95)$.
\bigskip

\noindent {\stampatello Corollary 1.} {\it Let $\theta\in (0,1/2)$ be fixed. If $f:\N \rightarrow \R$ has level $\lambda\in\big(0,\max\{(1+\theta)/2, 48/95\}\big)$, then there exists $\varepsilon_0=\varepsilon_0(\theta,\lambda)>0$ such that, $$
J_f(N,h)\ll_{\varepsilon_0} N^{1-\varepsilon_0}h^2, \qquad I_f(N,h)\ll_{\varepsilon_0} N^{1-\varepsilon_0}h^2,
$$
as $N\to \infty$, where $h=[N^{\theta}]$.}
\bigskip

\noindent
Note that $48/95=1/2+1/190>(1+\theta)/2$ if and only if $\theta<1/95$. Unlike [C1], we derive
Theorem 1 from a slight generalization concerning the {\it mixed} Selberg integral and 
the {\it mixed} symmetry integral of the sieve functions $f_1$ and $f_2$, namely
$$
J_{f_1,f_2}(N,h)\defineq \int_{N}^{2N}\prod_{c=1,2}\Big(\sum_{x<n\le x+h}f_c(n) - M_{f_c}(h)\Big)\,{\rm d}x,
$$
$$
I_{f_1,f_2}(N,h)\defineq \int_{N}^{2N} \prod_{c=1,2}\Big( \sum_{|n-x|\le h}\sgn(n-x)f_c(n)\Big) \,{\rm d}x. 
$$
\par
\noindent
where as before we set
$$
M_{f_c}(h)\defineq h\sum_{d\le Q_c}{{g_c(d)}\over d}, 
\quad 
(c=1,2),
$$ 
provided that $g_c$ and $Q_c$ are the Eratosthenes transform and the range of $f_c$, respectively.
\bigskip

\noindent {\stampatello Theorem 2.} {\it Fix  a small $\delta>0$. If for each $c=1,2$ the real sieve function $f_c$ has level
$\lambda_c\in(1/2,1)$ with $\lambda_1\ge \lambda_2$, then
$$
J_{f_1,f_2}(N,h)\EssBdd Nh+N^{\delta}Q_1^{53/48}Q_2^{7/8}h^2+N^{1-2\delta/3}h^2+Q_1h^2,
$$
$$ 
I_{f_1,f_2}(N,h)\EssBdd Nh+N^{\delta}Q_1^{53/48}Q_2^{7/8}h^2+N^{1-2\delta/3}h^2,
$$
where $Q_c=[N^{\lambda_c}]$ and $h=[N^{\theta}]$ with $\theta\in (0,1/2)$.}
\bigskip

It is plain that Theorem 1 follows immediately by taking $\lambda_1=\lambda_2=\lambda$, $Q_1=Q_2$, and $f_1=f_2=f$
in Theorem 2. On the other side, since Theorem 1.1 and Corollary 1.1 of [C1]
can be easily extended to 
				%PAGE 3
$J_{f_1,f_2}(N,h)$ and $I_{f_1,f_2}(N,h)$, then we 
can combine such a generalization with Theorem 2  to get the following immediate consequence.

\bigskip

\noindent {\stampatello Corollary 2.} {\it Let $\theta\in (0,1/2)$ be fixed and let $\lambda_1\ge \lambda_2>0$ be such that $\lambda_1+\lambda_2<1$ or  $53\lambda_1+42\lambda_2<48$. If for each $c=1,2$ the real sieve function $f_c$ has level
$\lambda_c$, then there exists $\varepsilon_0=\varepsilon_0(\theta,\lambda_1,\lambda_2)>0$ such that
$$
J_{f_1,f_2}(N,h)\ll_{\varepsilon_0} N^{1-\varepsilon_0}h^2, \qquad I_{f_1,f_2}(N,h)\ll_{\varepsilon_0} N^{1-\varepsilon_0}h^2,
$$
as $N\to \infty$, where $h=[N^{\theta}]$.
}
\bigskip

Summarizing, we need to prove only Theorem 2 and this is accomplished in \S4. To this end,
we premise a short section on some further notation and basic formul\ae, where we introduce the crucial auxiliary function ${\cal R}(a)$ in terms 
of the first Bernoulli periodic function. In \S3 we 
give the necessary lemmata for Theorem 2. The first lemma is a famous theorem of Duke, Friedlander and Iwaniec and it is the novelty of the present approach to
estimating ${\cal R}(a)$. Such an estimate is the theme of the second lemma. The link to $J_{f_1,f_2}(N,h)$ and $I_{f_1,f_2}(N,h)$ is provided by
the {\it correlations} of the sieve functions $f_1,f_2$ for which 
the third and last lemma gives a 
formula, 
with an error term taken under control by the new bound of ${\cal R}(a)$. We conclude the paper with a section of further comments and with an appendix including the proof of the Fourier expansion of the  first Bernoulli periodic function on the rational numbers. 

\bigskip

\par
\noindent {\bf 2. Some further notation and recurrent properties}
\bigskip
As usual in number theory, $(m,n)$ denotes the greatest common divisor of 
integers $m$ and $n$. 
Although 
$(x,y)$ denotes also the pair with coordinates $x,y$ or the open interval with real endpoints $x,y$, the meaning will be evident from the context. For the same sake of brevity, we use to write  $n\equiv a\ (q)$
instead of $n\equiv a\ (\bmod\, q)$. Moreover, we set 
$e(\alpha)\defineq e^{2\pi i\alpha}$ $\forall \alpha\in \R$ 
and $e_{q}(n)\defineq e(n/q)$ $\forall (n,q)\in \Z\times\N$. 

\noindent
The {\it distance} of  $\alpha \in \R$ from the integers is 
$\displaystyle{\Vert \alpha \Vert \defineq \min_{n\in \Z}|\alpha -n|}$ and
$\{ \alpha \}\defineq \alpha - [\alpha]$ 
is its {\it fractional part}. For the main variable $N$ we set $L
\defineq 
\log N$.

Without further references, throughout the paper we will appeal to the well-known inequalities
$$
\sum_{V_1<v\le V_2}e(v\alpha)\ll \min\Big(V_2-V_1,{1\over\Vert\alpha\Vert}\Big),\qquad
\sum_{d|t}1\ll t^{\varepsilon}\ (\forall t\in\N, \forall \varepsilon>0).
$$
Let us recall that the {\it first Bernoulli periodic function} is defined as

$$
{\cal B}_1(\alpha)\defineq\cases{\{ \alpha \} - 1/2\ &\hbox{if}\ $\alpha \in \R \backslash \Z$,\cr 0\ &\hbox{otherwise},\cr}
$$
whose (finite) Fourier expansion on the rational numbers is given by  (see 
the 
Appendix for the proof)
$$
{\cal B}_1\Big( {n\over q}\Big)=-{1\over q}\sum_{j\le {q\over 2}}\cot {{\pi j}\over q} \sin {{2\pi jn}\over q} \qquad \forall (n,q)\in \Z\times\N\setminus\{1\}.
\eqno(1)
$$
See that
$\cot(\pi/2+k\pi)=0$, $\forall k\in\Z$. 
Such expansion is interpreted as ${\cal B}_1=0$ when $q=1$.
It is easy to see that, for any $\alpha\in (0,+\infty)\backslash \N$, $d,q\in\N$ one 
has
$$
\#\{m\in(\alpha,2\alpha]:\ m\equiv d\, (q)\}=
{{[2\alpha]-[\alpha]}\over q}+\cases{\displaystyle{{\cal B}_1\Big({{[\alpha]-d}\over q}\Big)-{\cal B}_1\Big({{[2\alpha]-d}\over q}\Big)}&\hbox{if}\ $q\not \vert \, [c\alpha]-d$ for $c=1,2$,\cr\cr O(1)\ &\hbox{otherwise}.\cr}
\eqno(2)
$$

\noindent
Given the functions $g_1, g_2$ supported in $[1,Q_1], [1,Q_2]$, respectively, for 
all $a\in \Z\setminus\{0\}$ we set
$$
{\cal R}(a)\defineq\sum_{{\ell|a}\atop {}}\sum_{\qone  \sim {Q_1\over {\ell}}}g_1(\ell \qone  )\sum_{{\qtwo \sim {Q_2\over {\ell}}}\atop {(\qone  ,\qtwo )=1}}g_2(\ell \qtwo )\sum_{c=1,2}
(-1)^{c+1}\Big( {\cal B}_1\Big({{[\alpha_c]+\qoneinv  b}\over \qtwo}\Big) +{\cal B}_1\Big({{[\alpha_c]-\qoneinv  b}\over \qtwo}\Big) \Big),
$$
\par				%PAGE 4
\noindent
where 
$x\sim X$ means that $x\in(X,2X]\cap\N$, the integer $\qoneinv\in[1,\qtwo]$ is defined by $\qoneinv  \qone  \equiv 1\ (\bmod\, \qtwo )$
when $(\qone  ,\qtwo )=1$, and we set $b\defineq |a|/\ell$, $\alpha_c\defineq cN/\ell\qone$.
We explicitly remark that
${\cal R}(a)$ depends also on $g_1,g_2$ and $N$.
\noindent
Note that  $Q_c\ge |a|$ ensures $(Q_c/\ell,2Q_c/\ell]\cap\N\not=\emptyset$ for every $\ell|a$. On the other side, 
$Q_1\ll |a|$ yields 
${\cal R}(a)\EssBdd |a|Q_2$, which in turn becomes ${\cal R}(a)\EssBdd N^{1-\delta}$ when we assume also $Q_2=o(N^{1-\delta}/|a|)$
(the same property holds by interchanging the roles of $Q_1$ and $Q_2$). 

\bigskip

\par
\noindent {\bf 3. Lemmata}
\bigskip
The first lemma comes in a straightforward way from Theorem 2 of [DFI]. 
\bigskip

\noindent {\stampatello Lemma 1}. {\it Let $N, Q_1,Q_2\in \N$ and $k\in \Z\setminus\{0\}$ such that 
$Q_1, Q_2\le N$ and $k\ll Q_1Q_2$, as $Q_1, Q_2\to \infty$. If $g_1,g_2: \N \rightarrow \R$ are essentially bounded and supported in 
$[Q_1,2Q_1]$ and $[Q_2,2Q_2]$, 
respectively, then} 
$$
\sum_{\qone  \sim Q_1}g_1(\qone)\sum_{{\qtwo \sim Q_2}\atop {(\qone,\qtwo)=1}}g_2(\qtwo)e_{\qtwo}\!(k\, \qoneinv)
\EssBdd (Q_1Q_2)^{7\over 8}(Q_1+Q_2)^{{11}\over {48}}. 
$$

The first part of next lemma is basically Lemma B of [C2], 
reformulated for a pair of essentially bounded functions with support in a bounded interval. Here we take this opportunity to give a more detailed proof. The second part is where we apply the previous lemma and it constitutes the novelty of the 
method, 
in comparison with [C1] and [C2].

\bigskip

\noindent {\stampatello Lemma 2}. {\it Fix  a sufficiently small $\delta>0$. Let  $N, Q_1,Q_2\in \N$ and $a\in \Z\setminus\{0\}$ such that $Q_1Q_2\gg N^{1-2\delta/3}$, $Q_2\ll Q_1=o(N^{1-\delta})$ and $|a|=o(N)$,  as $N\to \infty$. If $g_1,g_2: \N \rightarrow \R$ are essentially bounded and supported in 
$[Q_1,2Q_1]$ and $[Q_2,2Q_2]$, 
respectively, then for every $\varepsilon>0$ one has
$$
{\cal R}(a) = {2\over {\pi}}\sum_{\ell|a}\sum_{\qone  \sim Q_1/\ell}g_1(\ell \qone  )\sum_{{\qtwo \sim Q_2/\ell}\atop {(\qone  ,\qtwo )=1}}g_2(\ell \qtwo )\sum_{j\le J}{\Delta_j\over j}
+ O_{\varepsilon}( N^{1-\delta+\varepsilon}),
\leqno{I)}
$$
\par
\noindent
where $J=J(\ell,\qone,\qtwo,N,\delta)\defineq [\ell \qone\qtwo N^{\delta-1}]$ and 
$$
\Delta_j=\Delta_j(\ell,\qone,\qtwo,a,N)\defineq\Big(\sin {{2\pi [2N/\ell \qone]j}\over \qtwo}-\sin {{2\pi [N/\ell \qone]j}\over \qtwo}\Big)\cos {{2\pi\, \qoneinv  |a|j}\over \ell \qtwo}.
$$
Also,  
$$
{\cal R}(a)\ll_{\varepsilon} N^{\varepsilon}(N^{\delta}Q_1^{53/48}Q_2^{7/8} + N^{1-\delta}). 
\leqno{II)}
$$
}

\noindent {\stampatello proof}.  I) For every $\ell|a$, let us set $b\defineq|a|/\ell\in\N$ and note that  $J=o(Q_2/\ell)$, while
$M\defineq Q_1Q_2N^{\delta-1}L^{-1}\geq 1$ from the hypothesis $Q_1Q_2\gg N^{1-2\delta/3}$. 
Moreover, since ${\cal B}_1\ll 1$, 
the contribution to ${\cal R}(a)$ from all $\ell|a$ such that
$\ell>M$ is  trivially
$$
\EssBdd Q_1Q_2\sum_{{\ell|a}\atop {\ell > M}}\,{{1}\over {\ell^2}}\, 
\EssBdd N^{1-\delta}. 
$$
Thus, let us consider the sum over $\ell\le M$ such that $\ell|a$. Together with $\qone  \sim Q_1/\ell$ and $\qtwo \sim Q_2/\ell$,  the condition $\ell\le M$ yields $J=[\ell \qone  \qtwo N^{\delta-1}]\to \infty$, as $N\to\infty$. 
By using formula (1) and  the identity 
$\sin(x-w)-\sin(y-w)+\sin(x+w)-\sin(y+w)=2(\sin x-\sin y)\cos w$, it is readily seen that 

$$
{\cal R}(a)=O_{\varepsilon}\Big( N^{1-\delta+\varepsilon}\Big)+2\sum_{{\ell|a}\atop{\ell \le M}}\sum_{\qone  \sim Q_1/\ell}g_1(\ell \qone  )\sum_{{\qtwo \sim Q_2/\ell}\atop {(\qone  ,\qtwo )=1}}{{g_2(\ell \qtwo )}\over \qtwo}\sum_{j\le \qtwo /2}
\Delta_j\cot {{\pi j}\over \qtwo}.\eqno(3)
$$
				%PAGE 5
Let us split the last sum as
$$
\sum_{j\le \qtwo /2}
\Delta_j\cot {{\pi j}\over \qtwo}=\sum_{j\le J}\Delta_j\cot {{\pi j}\over \qtwo}+\sum_{J<j\le \qtwo /2}\Delta_j\cot {{\pi j}\over \qtwo}=
{\cal D}_1+{\cal D}_2,\ \hbox{say},
$$
and first evaluate the contribution to (3) from ${\cal D}_2$. To this end, note that
$$
\eqalign{
\Delta_j=&{1\over{2i}}\cos\Big({2\pi\qoneinv bj\over \qtwo}\Big)\sum_{c=1,2}(-1)^c\Big(
e_{\qtwo}\!(j[\alpha_c])-e_{\qtwo}\!(-j[\alpha_c])\Big)\cr
=&{1\over {4i}}\Big(
{\cal E}_{2,j}-{\cal E}_{1,j}+\overline{\cal E}_{1,j}-\overline{\cal E}_{2,j}\Big),
\cr
}
$$
where for $c=1,2$ we set $\alpha_c\defineq cN/\ell \qone$ and
$$
{\cal E}_{c,j}={\cal E}_{c,j}(\ell,\qone,\qtwo,b,N)\defineq e_{\qtwo}\!\big(j([\alpha_c]+\qoneinv  b)\big)+e_{\qtwo}\!\big(j([\alpha_c]-\qoneinv  b)\big).
$$
Since $\cot(\pi/2)=0$,
by partial summation we can write
$$
\eqalign{
{\cal D}_2
&\ll \Big|\int_{J}^{\qtwo/2}\Big(\sum_{J<j\le v}\Delta_j\Big)\Big({{\rm d}\over {{\rm d}v}}\cot {{\pi v}\over \qtwo}\Big){\rm d}v\Big|\cr
&\ll \qtwo\int_{J}^{\qtwo/2}{1\over {v^2}}\Big|\sum_{J<j\le v}\Big(
{\cal E}_{2,j}-{\cal E}_{1,j}+\overline{\cal E}_{1,j}-\overline{\cal E}_{2,j}\Big)\Big|{\rm d}v.\cr
}
$$
The contribution from $\qtwo|\displaystyle{\qone[\alpha_c]\pm b}$ is trivially
$\displaystyle{\ll \qtwo\Big|\int_{J}^{\qtwo/2}{{{\rm d}v}\over {v}}\Big|\EssBdd \qtwo}$, which in turn contributes
to (3) as
$$
\eqalign{
&\EssBdd \sum_{\ell|a}\sum_{\qone  \sim Q_1/\ell}\#\{\qtwo \sim Q_2/\ell:\, (\qtwo,\qone)=1
\enspace \hbox{and} \enspace 
\displaystyle{\qtwo|\qone[\alpha_c]\pm b}\}\cr
&
\EssBdd \sum_{\ell|a}\Big( \sum_{{\qone  \sim Q_1/\ell}\atop {\qone  \not\,|\,b}}\sum_{
{(\qtwo ,\qone  )=1}\atop{\qtwo |\qone[\alpha_c]\pm b}}1+{Q_2\over {\ell}}\sum_{\qone  |b}1\Big)
\EssBdd (Q_1+Q_2)\sum_{\ell|a}{1\over\ell}\EssBdd N^{1-\delta},
\cr
}
$$
where we have 
taken into account: 
$\displaystyle{\qone[\alpha_c]\pm b\leq (cN\pm|a|)/\ell\ll N/\ell}$, $0<|a|\ll N$ 
 and  
$Q_2\ll Q_1=o(N^{1-\delta})$.

\noindent
On the other side, the contribution to ${\cal D}_2$ from $\displaystyle{\qone[\alpha_c]\pm b\not\equiv 0}$ (mod $\qtwo$) amounts to
$$
\ll
{\qtwo \over J}\sum_{c=1,2}
\Big( \Big \Vert {{[\alpha_c]+\qoneinv  b}\over \qtwo}\Big \Vert^{-1} +\Big \Vert {{[\alpha_c]-\qoneinv  b}\over \qtwo}\Big \Vert^{-1}\Big).
$$
Now, observe that $\displaystyle{\qone[\alpha_c]\pm b\not\equiv 0}$ (mod $\qtwo$), with $(\qone, \qtwo)=1$, yields
$\qone[\alpha_c]\pm b\equiv r\qone$ (mod $\qtwo$) for some $r$ such that  $1\le |r|\le \qtwo/2$, i.e. 
$$
\Big\Vert {{[\alpha_c]\pm\qoneinv  b}\over \qtwo}\Big \Vert={|r|\over \qtwo}.
$$
				%PAGE 6
Therefore, since $J=[\ell \qone  \qtwo N^{\delta-1}]$, the contribution to (3) from $\displaystyle{\qone[N/\ell \qone]+b\not\equiv 0}$ (mod $\qtwo$) through ${\cal D}_2$ is
$$
\eqalign{
\EssBdd &{{N^{1-\delta}}\over {Q_1Q_2}}\sum_{\ell|a}\ell \sum_{\qone  \sim {Q_1\over\ell}}
\sum_{{{\qtwo \sim Q_2/\ell}\atop {(\qtwo ,\qone  )=1}}\atop
{\qone[N/\ell \qone]+b\not \equiv 0\, (\qtwo)}}\Big \Vert {{[N/\ell \qone]+\qoneinv  b}\over \qtwo}\Big \Vert^{-1}\cr
\EssBdd &{{N^{1-\delta}}\over {Q_1Q_2}}\sum_{\ell|a}\ell \sum_{\qtwo \sim {Q_2\over\ell}}\qtwo
\sum_{1\le |r|\le {\qtwo\over 2}}{1\over |r|}\sum_{{{\qone  \sim Q_1/\ell}\atop {(\qone  ,\qtwo )=1}}\atop
{\qone[N/\ell \qone]+b\equiv r\qone (\qtwo)}}1 
\cr
\EssBdd &{{N^{1-\delta}}\over {Q_1Q_2}}\sum_{\ell|a}\ell \sum_{\qone  \sim {Q_1\over\ell}}
\sum_{1\le |r|\le {Q_2\over \ell}}{1\over |r|}
\sum_{{{\qtwo \sim Q_2/\ell}\atop {\qtwo| ([N/\ell \qone]-r)\qone+b}}
}\qtwo
\cr
\EssBdd &{{N^{1-\delta}Q_2}\over {Q_1}}\sum_{\ell|a}{1\over\ell} \sum_{\qone  \sim {Q_1\over\ell}\atop \qone\!|b}
\sum_{1\le |r|\le {Q_2\over \ell}\atop ([N/\ell \qone]-r)\qone+b=0}{1\over |r|}\cr
&+{{N^{1-\delta}}\over {Q_1}}\sum_{\ell|a}{1\over\ell^\varepsilon}\sum_{\qone  \sim {Q_1\over\ell}}
\sum_{1\le |r|\le {Q_2\over \ell}\atop ([N/\ell \qone]-r)\qone+b\not=0}{1\over |r|}
\cr
\EssBdd &{{N^{1-\delta}Q_2}\over {Q_1}}+N^{1-\delta}\EssBdd  N^{1-\delta}.\cr
}
$$
The same bound holds in the other three cases 
$\displaystyle{\qone[2N/\ell \qone]\pm b\not\equiv 0\, (\bmod \; \qtwo)}$,
$\displaystyle{\qone[N/\ell \qone]- b\not\equiv 0\, (\bmod \; \qtwo)}$.

Now, we turn our attention to ${\cal D}_1$. By using  the expansion of the cotangent function in power series, for a fixed $K>1$ we can write 
$$
{\cal D}_1= {\qtwo \over {\pi}}\sum_{j\le J}{\Delta_j\over {j}}+\sum_{n=1}^{K-1}{a_n\over\qtwo\!^n}\sum_{j\le J}j^n\Delta_j
+ O\Big({J\Big({J\over {\qtwo}}\Big)^{K}}\Big) ,
$$
where the coefficients $a_n\ll 1$ are given in terms of the Bernoulli numbers  (see  [MV], Appendix B, exercise 11 and formula (B.20)). Note that the first sum gives the main term of the stated formula for ${\cal R}(a)$, whereas 
the $O$-term is $\ll q_2(Q_1/N^{1-\delta})^{K+1}\ll q_2$ from the hypothesis $Q_1=o(N^{1-\delta})$. In order to see that also the sum over $n$ contributes to (3) as a remainder term, we first apply
partial summation to write
$$
\sum_{n=1}^{K-1}{a_n\over\qtwo\!^n}\sum_{j\le J}j^n \Delta_j 
\ll \Big(\sum_{n=1}^{K-1}{J^n\over\qtwo\!^n}\Big)\max_{v\le J}\Big| \sum_{j\le v}\Delta_j\Big|.
$$
Then we observe that 
the argument previously used for ${\cal D}_2$ applies here, because it turns out that
$$
\sum_{n=1}^{K-1}{J^n\over\qtwo\!^n}\ll \sum_{n=1}^{K-1}\Big({Q_1\over {N^{1-\delta}}}\Big)^n\ll 1.
$$ 
The formula 
I) is completely proved.

In order to prove the inequality stated in II), from what we have seen in the proof of I), it is plain that we may confine to  consider only
$$
{\cal R}_c(a)\defineq\sum_{{\ell|a}\atop {\ell \le {M}}}\sum_{\qone  \sim {Q_1\over {\ell}}}g_1(\ell \qone  )\sum_{{\qtwo \sim {Q_2\over {\ell}}}\atop {(\qone  ,\qtwo )=1}}g_2(\ell \qtwo )
\sum_{j\le 4{{ML}\over \ell }}
{\Sigma_j(c)\over j},\quad (c=1,2),
$$
				%PAGE 7
where $M=Q_1Q_2N^{\delta-1}L^{-1}$ is defined as above, and
$$
\Sigma_j(c)\defineq\sin{{2\pi [\alpha_c]j}\over\qtwo}\cos {{2\pi \qoneinv  b j}\over \qtwo},\quad\hbox{with}\ \alpha_c={{cN}\over{\ell \qone}}
,\  b={{|a|}\over {\ell}}.
$$
Thus, we have to show that
${\cal R}_c(a)\EssBdd N^{\delta}Q_1^{53/48}Q_2^{7/8} + N^{1-\delta}$, $\forall c=1,2$.
To this end, we  write
$$
\eqalign{
\Sigma_j(c)=&\sin {{2\pi \alpha_cj}\over {\qtwo}}\cos {{2\pi \qoneinv  b j}\over \qtwo}\cr
&-(1-\cos {{2\pi \{\alpha_c\}j}\over \qtwo})
\sin {{2\pi \alpha_cj}\over {\qtwo}}\cos {{2\pi \qoneinv  b j}\over \qtwo}\cr
&-\sin {{2\pi \{\alpha_c\}j}\over \qtwo}\cos {{2\pi \alpha_cj}\over {\qtwo}}\cos {{2\pi \qoneinv  b j}\over \qtwo}\cr
=&\Sigma_j^{(0)}(c)-\Sigma_j^{(1)}(c)-\Sigma_j^{(2)}(c),\ \hbox{say.}
} 
$$
Accordingly we have ${\cal R}_c(a)={\cal R}_c^{(0)}(a)-{\cal R}_c^{(1)}(a)-{\cal R}_c^{(2)}(a)$ with 
$$
{\cal R}^{(\nu)}_c(a)\defineq\sum_{{\ell|a}\atop {\ell \le {M}}}\sum_{\qone  \sim {Q_1\over {\ell}}}g_1(\ell \qone  )\sum_{{\qtwo \sim {Q_2\over {\ell}}}\atop {(\qone  ,\qtwo )=1}}g_2(\ell \qtwo )
\sum_{j\le 4{{ML}\over \ell }}
{\Sigma_j^{(\nu)}(c)\over j},\ (\nu=0,1,2).
$$
By applying partial summation with respect to $\qone$ we see that
$$
\eqalign{
\sum_{\qone  \sim {Q_1\over {\ell}}}g_1(\ell \qone  )\sum_{{\qtwo \sim {Q_2\over {\ell}}}\atop {(\qone  ,\qtwo )=1}}&g_2(\ell \qtwo )
\Sigma_j^{(0)}(c)= 
\sum_{\qone  \sim {Q_1\over {\ell} }}g_1(\ell \qone  )\sum_{{\qtwo \sim {Q_2\over {\ell} }}\atop {(\qone  ,\qtwo )=1}}g_2(\ell \qtwo )
   \sin {{\pi cNj}\over {Q_1\qtwo}}\cos {{2\pi \qoneinv  bj}\over \qtwo} \cr
&+ {{2\pi c Nj}\over \ell}\int_{Q_1/\ell}^{2Q_1/\ell}\sum_{{Q_1\over {\ell} }<\qone\le v}g_1(\ell \qone  )\sum_{{\qtwo \sim {Q_2\over {\ell} }}\atop {(\qone  ,\qtwo )=1}}{{g_2(\ell \qtwo )}\over \qtwo}\cos {{2\pi \qoneinv  bj}\over \qtwo} \cos {{2\pi cNj}\over {\ell\qtwo v}}{{{\rm d}v}\over {v^2}}.\cr
}
$$
Thus, from Lemma 1 it follows that ${\cal R}^{(0)}_c(a)\EssBdd N^{\delta}Q_1^{53/48}Q_2^{7/8}$. 
Now let us prove ${\cal R}^{(1)}_c(a),{\cal R}^{(2)}_c(a)\EssBdd N^{1-\delta}$. To this end, 
by using  the expansion of the cosine function in power series, we fix an integer $K>1$ and write
$$
\eqalign{
\sum_{j\le 4ML/\ell}{\Sigma_j^{(1)}(c)\over j}=&\sum_{j\le 4ML/\ell}\Big(1-\cos {{2\pi \{\alpha_c\}j}\over \qtwo}\Big){\Sigma_j^{(0)}(c)\over j}\cr
=& 
\sum_{n=1}^{K-1}{b_n\over\qtwo\!^n}\sum_{j\le 4ML/\ell}j^{n-1}\Sigma_j^{(0)}(c)+ O\Big(\Big({{ML\over {\ell\qtwo}}\Big)^{K}}\Big)\cr
\ll& {\ell\over ML}\max_{v\le 4ML/\ell}\Big|\sum_{j\le v}\Sigma_j^{(0)}(c)\Big|\sum_{n=1}^{K-1}\Big({4ML\over\ell\qtwo}\Big)^n\cr
\ll& {\ell N^{1-\delta}\over Q_1Q_2}\max_{v\le {4Q_1Q_2\over \ell N^{1-\delta}}}\Big|\sum_{j\le v}
({\cal E}'_{c,j}-\overline{\cal E}'_{c,j})\Big|,\cr
}
$$
where we have set ${\cal E}'_{c,j}\defineq e_{\qtwo}\!\big(j(\alpha_c+\qoneinv  b)\big)+e_{\qtwo}\!\big(j(\alpha_c-\qoneinv  b)\big)$, and we have used 
$b_n\ll 1$, $ML/(\ell\qtwo)=o(1)$ (the latter following straightforwardly from the hypothesis $Q_1=o(N^{1-\delta})$). 
Now, it is easy to see that, according as $\qtwo|\alpha_c\pm\qoneinv  b$ or not, the same arguments adopted in the proof of I) to treat the exponential sums lead to
$$
{\cal R}^{(1)}_c(a)\EssBdd{N^{1-\delta}\over Q_1Q_2}\sum_{{\ell|a}\atop {\ell \le {M}}}\ell\sum_{\qone  \sim {Q_1\over {\ell}}}\sum_{{\qtwo \sim {Q_2\over {\ell}}}\atop {(\qone  ,\qtwo )=1}}
\max_{v\le {4Q_1Q_2\over \ell N^{1-\delta}}}\Big|\sum_{j\le v}
({\cal E}'_{c,j}-\overline{\cal E}'_{c,j})\Big|\EssBdd N^{1-\delta}.
$$
				%PAGE 8
In a completely similar way,  we conclude also that ${\cal R}^{(2)}_c(a)\EssBdd N^{1-\delta}$
after using the expansion of the sine function in power series and noticing that 
$$
4\cos {{2\pi \alpha_cj}\over {\qtwo}}\cos {{2\pi \qoneinv  b j}\over \qtwo}={\cal E}'_{c,j}+\overline{\cal E}'_{c,j}.
$$
The lemma is completely proved.\hfill \square
\bigskip

Now, we can state and prove the main lemma of the paper. It gives a fairly general asymptotic formula for 
the {\it correlation} 
$$
\Corr_{f_1,f_2}(a)\defineq \sum_{n\sim N}f_1(n)f_2(n-a) 
$$
of real sieve functions $f_1,f_2$. In particular, it provides 
a strong level ($>1/2$) for the autocorrelation $\Corr_{f}=\Corr_{f,f}$
of a real sieve function $f$. It is also worthwhile to remark that 
next Lemma 
applications to $J_{f_1,f_2}$ and $I_{f_1,f_2}$ (i.e., Theorem 2) improve 
the non-trivial bounds given in [C1].
\bigskip

\noindent 
{\stampatello Lemma 3.} {\it Fix  a sufficiently small $\delta>0$. Let  $N, Q_1,Q_2\in \N$ such that $Q_2\to \infty$  and $Q_2\ll Q_1\ll N^{1-\delta}$,  as $N\to \infty$. If $g_1,g_2: \N \rightarrow \R$ are essentially bounded and supported in $[1,Q_1]$ and $[1,Q_2]$, respectively, then for every $\varepsilon>0$ and uniformly $\forall a\in \Z\setminus\{0\}$, with $|a|=o(N)$, one has
$$
{{\Corr_{f_1,f_2}(a)+\Corr_{f_1,f_2}(-a)}\over 2}=N\sum_{\ell |a}{1\over {\ell}}
\doublesum_{{\qone  \le Q_1/\ell}\atop{{\qtwo  \le Q_2/\ell} \atop(\qone,\qtwo)=1}}
 {{g_1(\ell \qone  )g_2(\ell \qtwo )}\over \qone\qtwo}
 + O_{\varepsilon}\big(N^{\delta+\varepsilon}Q_1^{53/48}Q_2^{7/8}
 +N^{1-2\delta/3+\varepsilon}\big),
$$
where $f_1=g_1\ast \1$, $f_2=g_2\ast \1$. }
\bigskip

\noindent {\stampatello proof}. First, we observe that
$$
\eqalign{
\Corr_{f_1,f_2}(a)=&\sum_{n\sim N}f_1(n)
\sum_{{\qtwo|n-a}\atop {\qtwo\le Q_2}}g_2(\qtwo)
=\sum_{\qtwo\le Q_2}g_2(\qtwo)\sum_{{n\sim N}\atop {n\equiv a\, (\qtwo)}}\sum_{{\qone|n}\atop {\qone\le Q_1}}g_1(\qone)\cr
=&\sum_{\ell|a}\sum_{\qtwo\le Q_2}g_2(\qtwo)\sum_{{\qone\le Q_1}\atop{(\qone,\qtwo)=\ell}}g_1(\qone)
\sum_{
{{n\sim N}\atop {n\equiv 0\, (\qone\!)}}\atop {n\equiv a\, (\qtwo\!)}}1\cr
=&\sum_{\ell|a}\sum_{\qone  \le{Q_1\over {\ell} }}g_1(\ell \qone  )\sum_{{\qtwo \le {Q_2\over {\ell} }}\atop {(\qone  ,\qtwo )=1}}g_2(\ell \qtwo )
\sum_{{m\sim {N\over {\ell \qone  }}}\atop {m\equiv \qoneinv\! b\, (\qtwo)}}1,\cr
}
$$
where we set $b=|a|/\ell$ as before. Then, plainly  we can write
$$
{{\Corr_{f_1,f_2}(a)+\Corr_{f_1,f_2}(-a)}\over 2}={1\over 2}\sum_{\ell|a}\sum_{\qone  \le{Q_1\over {\ell} }}g_1(\ell \qone  )\sum_{{\qtwo \le {Q_2\over {\ell} }}\atop {(\qone  ,\qtwo )=1}}g_2(\ell \qtwo )\sum_{{m\sim {N\over {\ell \qone  }}}\atop {m\equiv \pm\qoneinv\! b\, (\qtwo)}}1.
$$
Now, let us set $Q_{c,k}\defineq 2^{-k-1}Q_c$ for all $k=0,\ldots,[\log_2Q_c]$, $c=1,2$,
and confine to the dyadic intervals 
$(Q_{c,k}/\ell,2Q_{c,k}/\ell]$, where we define ${\cal R}_k(a)$ analogously to ${\cal R}(a)$. Thus, 
by the formula (2) we get
$$
\eqalign{
{1\over 2}\sum_{\ell|a}\sum_{\qone  \sim{Q_{1,k}\over {\ell} }}g_1(\ell \qone  )\sum_{{\qtwo \sim {Q_{2,k}\over {\ell} }}\atop {(\qone  ,\qtwo )=1}}g_2(\ell \qtwo )\sum_{{m\sim {N\over {\ell \qone  }}}\atop {m\equiv \pm\qoneinv\! b\, (\qtwo)}}1
=&N\sum_{\ell |a}{1\over {\ell}}\doublesum_{\qc  \sim{Q_{c,k}\over {\ell} }\atop(\qone,\qtwo)=1}{{g_1(\ell \qone  )g_2(\ell \qtwo )}\over \qone\qtwo}
 + {{{\cal R}_k(a)}\over 2} 
\cr
&+O_{\varepsilon}\Big(N^{\varepsilon}\sum_{\ell|a}\sum_{\qone  \sim{Q_{1,k}\over {\ell}}}
\sum_{\qtwo \sim Q_{2,k}/\ell\atop \qtwo|\qone[\alpha_c]\pm b}1
\Big). \cr
}
$$
				%PAGE 9
Since $|a|=o(N)$ yields $\qone[cN/(\ell \qone  )]\pm b \neq 0$, the latter $O$-term is $\EssBdd Q_1\EssBdd N^{1-\delta}$.
Moreover, we can assume  that $Q_{1,k}Q_{2,k}\gg N^{1-2\delta/3}$, for otherwise trivially ${\cal R}_k(a)\EssBdd N^{1-2\delta/3}$. Hence, the conclusion follows from II) of Lemma 2. \hfill \square
\bigskip

\noindent {\stampatello Remark.} Since $f_1,f_2$ are  essentially bounded, then for any $a>0$ one has
$$
\Corr_{f_1,f_2}(-a)=\sum_{n\sim N}f_1(n)f_2(n+a)=\sum_{N+a<n\le 2N+a}f_2(n)f_1(n-a)=\Corr_{f_2,f_1}(a)
+O_{\varepsilon}(aN^{\varepsilon}),
$$
In particular, if $f_1=f_2 $, this implies that
$$
{{\Corr_{f_1,f_2}(a)+\Corr_{f_1,f_2}(-a)}\over 2}=\Corr_f(a)+O_{\varepsilon}\Big( N^{\varepsilon}|a|\Big).
$$
Therefore, from the previous lemma we obtain the following formula for the value attained at $a=o(N)$ by the autocorrelation of a sieve function 
$f=g\ast \1$  of range $Q\ll N^{1-\delta}$:
$$
\Corr_f(a) = N\sum_{\ell |a}{1\over {\ell}}
\doublesum_{{\qone,\qtwo  \sim Q/\ell}\atop{(\qone,\qtwo)=1}}
{{g(\ell \qone  )g(\ell \qtwo )}\over \qone\qtwo}
 + O_{\varepsilon}\Big(N^{\delta+\varepsilon}Q^{{95}\over {48}}+N^{1-2\delta/3+\varepsilon}\Big),
$$
that means level $\lambda=1/2+1/190$ for autocorrelations of $f$.
\bigskip

\par
\noindent {\bf 4. Proof of Theorem 2}
\bigskip

\noindent
Let us consider the symmetry integral first, and write
$$
\eqalign{
I_{f_1,f_2}(N,h)=&\doublesum_{N-h<n,m\le 2N+h\atop 0\le |n-m|\le 2h}f_1(n)f_2(m)\integrale_{{N\le x\le 2N}\atop {{|x-n|\le h}\atop {|x-m|\le h}}}\sgn(x-n)\sgn(x-m){\rm d}x\cr
=&\doublesum_{N+h<n,m\le 2N-h\atop 0\le |n-m|\le 2h}f_1(n)f_2(m)\integrale_{{|t|\le h}\atop {|t+(n-m)|\le h}}\sgn(t)\sgn(t+(n-m)){\rm d}t + O_{\varepsilon}(N^{\varepsilon}h^3) \cr
=&\sum_{N<n\le 2N}f_1(n)\sum_{0\le |a|\le 2h}f_2(n-a)\integrale_{{|t|\le h}\atop {|t+a|\le h}}\sgn(t)\sgn(t+a){\rm d}t + O_{\varepsilon}(N^{\varepsilon}h^3).
}
$$
\par
\noindent
To simplify our exposition, somewhere
the symbol (T) within some of the next formul\ae\ will warn the reader 
of some {\it tails}, i.e. terms being $\EssBdd h^3\EssBdd Nh$, that are discarded to abbreviate the formul\ae\ themselves. Thus, the above equation becomes
$$
I_{f_1,f_2}(N,h)\buildrel{\hbox{\piccolissimo (T)}}\over{\sim}
\sum_{0\le |a|\le 2h}W(a)\Corr_{f_1,f_2}(a),\quad 
\hbox{
with}\quad
W(a)\defineq \integrale_{{|t|\le h}\atop {|t-a|\le h}}\sgn(t)\sgn(t-a){\rm d}t.
$$
Since $f_1, f_2$ are essentially bounded, then 
$$
W(0)\Corr_{f_1,f_2}(0)=2h\sum_{n\sim N}f_1(n)f_2(n)\EssBdd Nh.
$$
Moreover, note that $W$ is even and $W(a)\ll h$ uniformly for all $a$. Therefore, Lemma 3 implies that
$$
\eqalign{
\sum_{0< |a|\le 2h}W(a)\Corr_{f_1,f_2}(a)=&\sum_{0<a\le 2h}W(a)(\Corr_{f_1,f_2}(a)+\Corr_{f_1,f_2}(-a))\cr
=&
2N\sum_{0<a\le 2h}W(a)\sum_{\ell |a}{1\over {\ell}}
\doublesum_{{\qone  \le Q_1/\ell}\atop{{\qtwo  \le Q_2/\ell} \atop(\qone,\qtwo)=1}}{{g_1(\ell \qone  )g_2(\ell \qtwo )}\over \qone\qtwo}\cr
&+
O_{\varepsilon}\big(N^{\delta+\varepsilon}Q_1^{53/48}Q_2^{7/8}h^2
 +N^{1-2\delta/3+\varepsilon}h^2\big).
 }
$$
				%PAGE 10
Hence, the stated inequality for $I_{f_1,f_2}(N,h)$ follows from
$$
\sum_{0<a\le 2h}W(a)\sum_{\ell |a}{1\over {\ell}}\doublesum_{{\qone  \le Q_1/\ell}\atop{{\qtwo  \le Q_2/\ell} \atop(\qone,\qtwo)=1}}{{g_1(\ell \qone  )g_2(\ell \qtwo )}\over \qone\qtwo}
\EssBdd\sum_{\ell \le 2h}{1\over {\ell}}\Big|\sum_{0<b\le 2h/\ell}W(\ell b)\Big|\EssBdd h,
$$
where we have applied the property (see [C1], Lemma 2.4) 
$$
\sum_{0<a\le 2h\atop {a\equiv 0\, (\ell)}}W(a)=2\ell\Big\Vert {h\over \ell}\Big\Vert \ll h, \quad \hbox{for}\ 1\le\ell\le 2h. 
$$

Now, let us turn our attention to the Selberg integral $J_{f_1,f_2}(N,h)$. First we observe that for any $c=1,2$ 
$$
\eqalign{
\int_{N}^{2N}\sum_{x<n\le x+h}f_c(n)\, {\rm d}x =& \sum_{N<n\le 2N+h}f_c(n)\integrale_{{N\le x\le 2N}\atop {n-h\le x<n}}{\rm d}x \cr
=& \sum_{N+h<n<2N-h}f_c(n)\int_{n-h}^{n}{\rm d}x + O_{\varepsilon}( N^{\varepsilon}h^2) \cr
=& h\sum_{n\sim N}f_c(n) + O_{\varepsilon}( N^{\varepsilon}h^2). 
}
$$
Since 
$$
\sum_{n\sim N}f_c(n) = \sum_d g_c(d)\Big( \Big[ {{2N}\over d}\Big] -\Big[ {N\over d}\Big] \Big) 
= N{M_{f_c}(h)\over h} + O_{\varepsilon}( N^{\varepsilon}Q_c),
$$
then, by recalling that $Q_1\geq Q_2$ and $M_{f_c}(h)\EssBdd h$, we can write
$$
J_{f_1,f_2}(N,h)
=\int_{N}^{2N}\doublesum_{x<n,m\le x+h}f_1(n)f_2(m){\rm d}x - NM_{f_1}(h)M_{f_2}(h)+O_{\varepsilon}\big( N^{\varepsilon}(Q_1h^2+h^3)\big).
$$
Now, by arguing as we have done before for $I_{f_1,f_2}(N,h)$, it is easy to see that
$$
\int_{N}^{2N}\doublesum_{x<n,m\le x+h}f_1(n)f_2(m){\rm d}x=
\sum_{0\le |a|\le h}(h-|a|)
\sum_{n\sim N}f_1(n)f_2(n+a)+ O_{\varepsilon}(N^{\varepsilon}h^3),
$$
which yields
$$
J_{f_1,f_2}(N,h)\buildrel{\hbox{\piccolissimo (T)}}\over{\sim}\sum_{0< |a|\le h}(h-|a|)\Corr_{f_1,f_2}(-a) - NM_{f_1}(h)M_{f_2}(h)+O_{\varepsilon}\big( N^{\varepsilon}(Nh+Q_1h^2)\big),
$$
Note that $h-|a|$ is an even function of the variable $a$. Thus, the previous calculations and Lemma 3 apply again here to get
$$
\eqalign{
\sum_{0< |a|\le h}(h-|a|)\Corr_{f_1,f_2}(-a)=&
\sum_{0< a\le h}(h-a)\big(\Corr_{f_1,f_2}(-a)+\Corr_{f_1,f_2}(a)\big)\cr
=&2N
\sum_{\ell \le h}{1\over {\ell}}\sum_{b\le h/\ell}(h-\ell b)
\doublesum_{{\qone  \le Q_1/\ell}\atop{{\qtwo  \le Q_2/\ell} \atop(\qone,\qtwo)=1}}{{g_1(\ell \qone  )g_2(\ell \qtwo )}\over \qone\qtwo}\cr
&+
O_{\varepsilon}\big(N^{\delta+\varepsilon}Q_1^{53/48}Q_2^{7/8}h^2
 +N^{1-2\delta/3+\varepsilon}h^2\big).
 }
$$
				%PAGE 11
By using the formula
$$
\sum_{b\le h/\ell}(h-\ell b)={{h^2}\over {2\ell}}+O(h), \quad \forall \ell \le h,
$$
we can write
$$
\eqalign{
2N\sum_{\ell \le h}{1\over {\ell}}\sum_{b\le h/\ell}(h-\ell b)
\doublesum_{{\qone  \le Q_1/\ell}\atop{{\qtwo  \le Q_2/\ell} \atop(\qone,\qtwo)=1}}{{g_1(\ell \qone  )g_2(\ell \qtwo )}\over \qone\qtwo}
&=Nh^2\sum_{\ell \le h}{1\over {\ell^2}}
\doublesum_{{\qone  \le Q_1/\ell}\atop{{\qtwo  \le Q_2/\ell} \atop(\qone,\qtwo)=1}}{{g_1(\ell \qone  )g_2(\ell \qtwo )}\over \qone\qtwo}
+
O_{\varepsilon}(N^{1+\varepsilon} h)\cr
&\buildrel{\hbox{\piccolissimo (T)}}\over{\sim}
N h^2 \sum_{\ell=1}^{\infty}{1\over {\ell^2}}\doublesum_{{\qone  \le Q_1/\ell}\atop{{\qtwo  \le Q_2/\ell} \atop(\qone,\qtwo)=1}}{{g_1(\ell \qone  )g_2(\ell \qtwo )}\over \qone\qtwo}\cr
&=
N h^2 \sum_{\ell=1}^{\infty}\doublesum_{{\qone  \le Q_1}\atop{{\qtwo  \le Q_2} \atop(\qone,\qtwo)=\ell}}{{g_1(\qone  )g_2(\qtwo )}\over \qone\qtwo}\cr
&=N h^2 \doublesum_{\qone\le Q_1\atop \qtwo  \le Q_2}{{g_1(\qone)}\over \qone}{{g_2(\qtwo )}\over \qtwo}= NM_{f_1}(h)M_{f_2}(h). \cr
 }
$$
Hence, we conclude that
$$
\sum_{0< |a|\le h}(h-|a|)\Corr_{f_1,f_2}(-a) - NM_{f_1}(h)M_{f_2}(h)
\ll_{\varepsilon}N^{\delta+\varepsilon}Q_1^{53/48}Q_2^{7/8}h^2
 +N^{1-2\delta/3+\varepsilon}h^2+N^{1+\varepsilon} h .
$$
Theorem 2 is completely proved.\enspace \square
\bigskip

\par
\noindent {\bf 5. Further comments and remarks}
\bigskip

\noindent
1. Though analogous definitions and results can be easily formulated for complex valued sieve functions, here we stick to the real case for simplicity. 

\noindent
2. A famous example of sieve function is given by the truncated divisor sum $\Lambda_R$, which has been exploited by
Goldston, Pintz and Y{\i}ld{\i}r{\i}m for their recent results [GPY]. We refer the reader to [CL2]
for an application of our recent study about the distribution of sieve functions to the case of $\Lambda_R$. 

\noindent
3.
The key of the present approach is the treatment of the {\it error term} ${\cal R}_f(a)$ arising 
from the formula of the autocorrelation of a sieve function $f=g\ast\1$ written for any nonzero integer $a=o(N)$ as
$$
\Corr_f(a)
= \sum_{\ell |a}\doublesum_{(d,q)=1}g(\ell \qone  )g(\ell \qtwo ){1\over \qtwo}\Big( \Big[ {{2N}\over {\ell \qone  }}\Big] - \Big[ {N\over {\ell \qone  }}\Big]\Big) + {\cal R}_f(a).
$$
In [C1] such an error term is defined by using the orthogonality of the additive characters. 
The estimate of ${\cal R}_f(a)$ leads to
the gain $\Delta=1/48$ given in Theorem 1 and, as noted above, it is essentially due to   the non-trivial bound of the bilinear forms with Kloosterman fractions by Duke, Friedlander and Iwaniec (see Lemma 1). In this respect, the recent improvement obtained by Bettin and Chandee 
[BC] would yield $\Delta=1/20$. Moreover, we think that this result might lead to an improvement of ours in respect of the short length $h$ as well. We are going to show such a further achievement in a future paper.

\noindent
4. In the literature the {\it level} of distribution of an arithmetic function $f$ in the arithmetic progressions is usually meant to be a positive real number $\lambda_{AP}(f)$ such that 
$$
\sum_{q\le Q}\max_{(a,q)=1}\Big| \sum_{{n\le x}\atop {n\equiv a\, (q)}}f(n)-{1\over {\varphi(q)}}\sum_{{n\le x}\atop {(n,q)=1}}f(n)\Big| 
\ll_{\varepsilon} {x^{1-\varepsilon}},\qquad (\varphi(q)\defineq |\{ a\le q, (a,q)=1\}|),
$$
\par				%PAGE 12
\noindent
holds for $Q\le x^{\lambda_{AP}(f)-\delta}$ (where $\delta>0$ is sufficiently small and depends on $\varepsilon>0$). For example, the celebrated Bombieri-Vinogradov Theorem gives $\lambda_{AP}(\Lambda)={1\over 2}$ for the von-Mangoldt function $\Lambda$. 
\par
It is a classical and standard argument (see [E] for example) to deduce from 
the previous inequality 
an asymptotic formula for the autocorrelations of a sieve function $f$ given by
$$
f(n)=\sum_{{q|n}\atop {q\le Q}}g(q), \quad\hbox{with}\quad Q\ll N^{\lambda_{AP}(f)-\delta}.
$$
Somehow this justifies our definiton of {\it level} when we refer to a sieve function.
Unfortunately, it seems to be  very hard to reverse such a process, that is to say,
information on the {\it level} of distribution of $f$ in the arithmetic progressions seems to be much stronger than knowledge about the autocorrelation of $f$. From this point of view, our Lemma 3 provides 
a substitute for what is still lacking from the study 
of distribution of $f$ in the arithmetic progressions.

\bigskip

\par
\noindent {\bf 6. Appendix: the first Bernoulli function on the rational numbers}
\bigskip

Here we prove the formula (1), that equivalently we state as
$$
{\cal B}_1\Big( {n\over q}\Big) = {i\over {2q}}
\sum_{0<|j|\le q/2}\cot {{\pi j}\over q}e_q(jn), \qquad \forall (n,q)\in \Z\times\N\setminus\{1\}.
$$
\noindent
To this end, it suffices to establish the following equality for the Fourier coefficient of ${\cal B}_1$:
$$
c_{j,q}\defineq{1\over q}\sum_{0\le |r|\le q/2}{\cal B}_1\Big( {r\over q}\Big)e_q(-jr)=
{i\over {2q}}\cot {{\pi j}\over q}.
$$
\noindent
First, since ${\cal B}_1$ is odd, note that $c_{j,q}={i\over q}\Sigma$, where we set 
$$
\Sigma\defineq -{2\over q}\sum_{r\le [q/2]}r\sin {{2\pi jr}\over q} 
+ \sum_{r\le [q/2]}\sin {{2\pi jr}\over q}.
$$
\noindent
Note that, for $R=[q/2]$ by Abel's lemma (see  [MV], Appendix A, exercise 3) one has
$$
\sum_{r\le R}r\sin {{2\pi jr}\over q} = R\sum_{r\le R}\sin {{2\pi jr}\over q} 
- \sum_{r\le R-1}\Big( \sum_{\ell \le r}\sin {{2\pi j\ell}\over q} \Big).
$$
Therefore, by using the identities (see [GR], formul\ae\ n.1.342.1 and n.1.342.2), $\forall X\in\N, \forall \alpha \in \R \backslash \Z$,
$$
\sum_{r\le X}\sin(2\pi \alpha r)= \sin^2(\pi \alpha X) \cot(\pi \alpha) 
+ {{\sin(2\pi \alpha X)}\over 2},
$$
$$
\sum_{r\le X}\cos(2\pi \alpha r)
= {{\sin(2\pi \alpha X) \cot(\pi \alpha)}\over 2} 
- {{1 - \cos(2\pi \alpha X)}\over 2},
$$
we can write 
$$
\eqalign{
\Sigma 
=& \Big( -{{2R}\over q}+1\Big) \sum_{r\le R}\sin {{2\pi jr}\over q} 
+ {2\over q}\sum_{r\le R-1} \Big( \sum_{\ell \le r}\sin {{2\pi j\ell}\over q}\Big) \cr
=& {2\over q} \cot {{\pi j}\over q} \sum_{r\le R-1}\sin^2 {{\pi jr}\over q} 
+ \Big( 1-{{2R}\over q}+{1\over q}\Big) \sum_{r\le R}\sin {{2\pi jr}\over q} 
- {1\over q} \sin {{2\pi jR}\over q}\cr
=&
{R\over q} \cot {{\pi j}\over q} 
+ {1\over q}\Big( \cos {{2\pi jR}\over q} - 1\Big) \cot {{\pi j}\over q} 
- {1\over {2q}} \cot {{\pi j}\over q} \Big( \sin {{2\pi jR}\over q} \cot {{\pi j}\over q} 
 + \cos {{2\pi jR}\over q} - 1\Big)\cr
&
+ {{1+2\{q/2\}}\over q} \Big( \sin^2 {{\pi jR}\over q} \cot {{\pi j}\over q} 
 + {1\over 2}\sin {{2\pi jR}\over q} \Big) 
- {1\over q} \sin {{2\pi jR}\over q},\cr
}
$$
				%PAGE 13
where 
$$
2\{q/2\}=q-2[q/2]=q-2R=\cases{1\ &if $q$ is odd,\cr 0\ &\hbox{otherwise}.\cr}
$$
Now, since from $2\pi jR= \pi j q- 2\pi j\{q/2\}$ it follows that

$$
\eqalign{
\cos {{2\pi jR}\over q} &= (-1)^j \cos {{2\pi j \{q/2\}}\over q},\cr
	\sin {{2\pi jR}\over q} &= (-1)^{j+1}2\Big\{{q\over 2}\Big\}\sin {{\pi j}\over q}, \cr
		\sin^2 {{\pi jR}\over q} &= {1\over 2} +{{(-1)^{j+1}}\over 2} \cos {{2\pi j \{q/2\}}\over q},\cr
		}
$$
then we get
$$
\eqalign{
\Sigma
=& {1\over 2}\cot {{\pi j}\over q} 
- 2\Big\{{q\over 2}\Big\}{1\over {2q}}\cot {{\pi j}\over q} 
+ {1\over {2q}}\Big( (-1)^j \cos {{\pi j 2\{q/2\}}\over q} - 1\Big)\cot {{\pi j}\over q} \cr
&
+ 2\Big\{{q\over 2}\Big\}{{(-1)^j}\over {2q}}\sin {{\pi j}\over q} \cot^2 {{\pi j}\over q} 
+ {{1+2\{q/2\}}\over q} \Big( {1\over 2} 
 - {{(-1)^j}\over 2}\cos {{\pi j 2\{q/2\}}\over q} \Big) \cot {{\pi j}\over q}\cr
=& {1\over {2q}}\cot {{\pi j}\over q} 
\Big( q - 2\Big\{{q\over 2}\Big\} + 2\Big\{{q\over 2}\Big\}(-1)^j\cos {{\pi j}\over q}
+ 2\Big\{{q\over 2}\Big\} \Big( 1 - (-1)^j \cos {{\pi j 2\{q/2\}}\over q} \Big) \Big) \cr
=& {1\over 2}\cot {{\pi j}\over q}.
}
$$

\par
\centerline{\bf References}
\bigskip
\smallskip
\item{\bf [BC]} S. Bettin, V. Chandee, {\it Trilinear forms with Kloosterman fractions},
preprint at arXiv:1502.00769v1.
\smallskip
\item{\bf [C1]} G. Coppola, {\it On the Correlations, Selberg integral and symmetry of sieve functions in short intervals}, J. Combinatorics and Number Theory {\bf 2.2} (2010), Article 1, 91--105.
\smallskip
\item{\bf [C2]} G. Coppola, {\it On the Correlations, Selberg integral and symmetry of sieve functions in short intervals, II}, Int. J. Pure Appl. Math. {\bf 58.3} (2010), 281--298.
\smallskip
\item{\bf [CL1]} G. Coppola, M. Laporta, {\it Symmetry and short interval mean-squares}, (submitted), preprint available at 
arXiv:1312.5701v2.
\smallskip
\item{\bf [CL2]} G. Coppola, M. Laporta, {\it Sieve functions in arithmetic bands}, (submitted), preprint available at \par
arXiv:1503.07502v2. 
\smallskip
\item{\bf [DFI]} W. Duke, J. Friedlander, H. Iwaniec, {\it Bilinear forms with Kloosterman fractions}, Invent. Math. {\bf 128} (1997), no. 1, 23--43.
\smallskip
\item{\bf [E]}
P.D.T.A. Elliott, {\it On the Correlation of Multiplicative and the Sum of Additive Arithmetic Functions}, Mem. Amer. Math. Soc., {\bf 112} (1994), no. 538.
\smallskip
\item{\bf [GPY]} D.A. Goldston, J. Pintz, C. Yildirim, {\it Primes in tuples. I},  Ann. of Math. (2) {\bf 170} (2009), no. 2, 819--862. 
\smallskip
\item{\bf [GR]}  I. S. Gradshteyn, I. M. Ryzhik, {\it Tables of integrals, series, and products}, Fifth Edition, Academic Press, 1994.  
\smallskip
\item{\bf [MV]}
H.L. Montgomery, R.C. Vaughan, {\it Multiplicative Number Theory. I. Classical Theory}, Studies in Advanced Math., vol. 97, Cambridge University Press, Cambridge, 2007.
\bigskip

\par
\leftline{\tt Giovanni Coppola \hfill Maurizio Laporta}
\leftline{\tt Universit\`a degli Studi di Napoli \hfill Universit\`a degli Studi di Napoli}
\leftline{\tt Home address \negthinspace : \negthinspace Via Partenio \negthinspace 12 \negthinspace - \hfill Dipartimento di Matematica e Appl.}
\leftline{\tt - 83100, Avellino(AV), ITALY \hfill Compl.Monte S.Angelo}
\leftline{\tt e-page : $\! \! \! \! \! \!$ www.giovannicoppola.name \hfill Via Cinthia - 80126, Napoli, ITALY}
\leftline{\tt e-mail : $\! \! \! \! \! \!$ giovanni.coppola@unina.it \hfill e-mail : mlaporta@unina.it}

\bye